\documentclass[11pt,a4paper]{article}
\usepackage{amsmath, amssymb}
\usepackage{latexsym, color}
\usepackage{epsfig}

\numberwithin{equation}{section}

\newtheorem{thm}{Theorem}[section]

\newtheorem{lem}[thm]{Lemma}
\newtheorem{prop}[thm]{Proposition}

\newcommand{\qed}{\hfill \ensuremath{\square}}
\newcommand{\ds}{\displaystyle}
\newcommand{\pf}{\noindent {\sl Proof}. \ }
\newcommand{\p}{\partial}

\newcommand{\norm}[1]{\| #1 \|}

\newcommand{\eqnref}[1]{(\ref {#1})}

\newcommand{\Rbb}{\mathbb{R}}

\newcommand{\Ncal}{\mathcal{N}}

\def\Ba{{\bf a}}

\def\Bc{{\bf c}}

\def\Be{{\bf e}}

\def\Bo{{\bf o}}
\def\Bp{{\bf p}}

\def\Bx{{\bf x}}
\def\By{{\bf y}}
\def\Bz{{\bf z}}

\newcommand{\Ga}{\alpha}
\newcommand{\Gb}{\beta}
\newcommand{\Gd}{\delta}
\newcommand{\Ge}{\epsilon}

\newcommand{\GD}{\Delta}

\newcommand{\GO}{\Omega}

\newcommand{\beq}{\begin{equation}}
\newcommand{\eeq}{\end{equation}}
\newcommand{\ol}{\overline}

\begin{document}

\title{Quantitative estimates of the field excited by an emitter
in a narrow region between two circular inclusions\thanks{\footnotesize This work was
supported by NRF grants No. 2015R1D1A1A01059212, 2016R1A2B4011304 and 2017R1A4A1014735, and by Hankuk University of Foreign Studies Research Fund of 2018.}}

\author{Hyeonbae Kang\thanks{\footnotesize Department of Mathematics and Institute of Applied Mathematics, Inha University, Incheon 22212, S. Korea (hbkang@inha.ac.kr).} \and KiHyun Yun\thanks{\footnotesize Department of Mathematics, Hankuk University of Foreign Studies, Yongin-si, Gyeonggi-do 17035, S. Korea (kihyun.yun@gmail.com).}}

\maketitle

\begin{abstract}
A field excited by an emitter can be enhanced due to presence of closely located inclusions. In this paper we consider such field enhancement when inclusions are disks of the same radii, and the emitter is of dipole type and located in the narrow region between two inclusions. We derive quantitatively precise estimates of the field enhancement in the narrow region. The estimates reveal that the field is enhanced by a factor of $\Ge^{-1/2}$ in most area, where $\Ge$ is the distance between two inclusions.  This factor is the same as that of gradient blow-up when there is a smooth back-ground field, not a field excited by an emitter. The method of deriving estimates shows clearly that enhancement is due to potential gap between two inclusions.
\end{abstract}



\section{Introduction and statements of results}

The purpose of this paper is to derive precise estimates for enhancement of fields excited by a dipole-type emitter in presence of closely located inclusions. In this paper we deal with the case when inclusions are of circular shapes as typical examples of domains with smooth boundaries to see how much the field is enhanced. This work is a continuation of its companion paper \cite{KY2} where the field enhancement is estimated when the inclusions are of a bow-tie shape with corners being separated by a small distance.

The problem of this paper and its companion paper is closely related to the study of enhancement of the smooth back-ground field, typically a uniform field, in the presence of  closely located inclusions. Such a study is motivated by the effective medium theory of densely packed perfect conductors \cite{FK-CPAM-73, Keller-JAP-63} and analysis of stress in composites with stiff inclusions \cite{bab, keller}. There has been significant progress on this subject in last two decades or so, for that we refer to \cite{KY2} and references therein.

In this paper we deal with the case when the field is excited by an emitter and enhanced by closely located inclusions, motivated by the study of nano-antennas (see, e.g., \cite{PBFLN}). The problem of this paper can be described in terms of the following mathematical model:
\beq\label{main_equation}
\begin{cases}
\GD u = {\bf a } \cdot \nabla  \Gd_{\Bp}  \quad&\mbox{in } \Rbb^2 \setminus \overline {(D_{1} \cup D_{2})}, \\
\ds u =c_j \quad&\mbox{on }\p D_{j}, \ j=1,2,  \\
\ds \int_{\p D_{j} } \p_{\nu} u \, ds  =0, ~&j=1,2,\\
\ds u (\Bx) = O\left( |\Bx|^{-1}\right)~&\mbox{as }|\Bx|\rightarrow\infty.
\end{cases}
\eeq
Here, $D_1$ and $D_2$ are bounded planar domains representing two inclusions whose distance is small, say $\Ge$:
$$
\Ge = \mbox{dist} (D_1, D_2).
$$
The second line in \eqnref{main_equation} requires that the solution $u$ attains constant values $c_j$ on $\p D_j$, which indicates that the inclusions $D_j$ are perfect conductors, namely, their conductivities are infinite. In general, the constant values are different, that is, $c_1 \neq c_2$, and this potential gap induces gradient blow-up or field enhancement. It is worthwhile to emphasize that $c_j$ are not prescribed, but to be determined by the problem. In particular, they depend on $\Ge$. The quantity ${\bf a } \cdot \nabla  \Gd_{\Bp}$ in the first line represents the emitter of the dipole type so that the unit vector $\Ba$ is the direction of the dipole and $\Bp$ its location. We assume that $\Bp$ is located in the narrow region in between two inclusions. Further, $\nu$ is the unit normal vector on $\p D_1 \cup \p D_2$ pointing inward to $D_1 \cup D_2$.

Let $\Ncal_\Bp(\Bx)$ be the fundamental solution to the Laplacian in two dimensions, namely,
\beq\label{Ncal}
\Ncal_\Bp(\Bx):= \frac{1}{2\pi} \log |\Bx-\Bp|.
\eeq
Then, $\GD \Ncal_\Bp(\Bx)= \Gd_{\Bp}(\Bx)$. Thus, in absence of inclusions, the solution to $\GD u = {\bf a } \cdot \nabla  \Gd_{\Bp}$ is given by
\beq\label{uzero}
\Ba \cdot \nabla \Ncal_\Bp(\Bx) = \frac{1}{2\pi} \frac{\Ba \cdot (\Bx-\Bp)}{|\Bx-\Bp|^2},
\eeq
and hence its gradient field has singularity at $\Bp$ of size $|\Bx-\Bp|^{-2}$. This singularity may be amplified by the interaction between two inclusions.

The objective of this paper is to estimate the gradient field in a narrow region in-between $D_1$ and $D_2$ when they are disks.
We assume that they are disks of same radii and their radii is of much larger scale than $\Ge$. Thus we suppose that their common radius is $1$. Then after translation and rotation we may assume that
\beq
D_j = B_{1}((-1)^j (1+\Ge/2),0), \quad j=1,2.
\eeq
Here and throughout this paper, $B_r (\Bc)$ denotes the open disk of radius $r$ centered at $\Bc$. If $\Bc$ is the origin $\Bo$, we simply write $B_r$. We assume that the emitter is located on the $x_2$-axis, i.e., $\Bp = (0, p)$ for some $p$ with $|p| \leq C$ for some $C$, say $C = 1/2$. We believe that similar results hold even if the radii are different even though analysis would be much more complicated technically, and it is not our intention here to pursue such a case.

The following are the main theorems of this paper in which (and throughout this paper) we employ commonly used symbols: the expression $\Ga \lesssim \Gb$ implies that there exists a positive constant $C$ independent of $\Ge$ (sufficiently small) such that $\Ga \leq C\Gb$, and $\Ga \simeq \Gb$ implies that both $\Ga \lesssim \Gb$ and $\Gb \lesssim \Ga$ hold.

\begin{thm}\label{thm_1st}
Let $u$ be the solution to \eqref{main_equation}. If $\Ba \neq  (0,1)$, then for any $M>0$ there exist positive constants $C_1$, $C_2$, $A_*$ and $\Ge_0$ depending on $\Ba$ and $M$ such that the following estimates hold for all ${\Bx}\in B_{1/2} \setminus \overline{ D_1 \cup D_2}$, for all $\Bp$ with $|\Bp| < M \sqrt \Ge$, and for all $\Ge< \Ge_0$:
\begin{itemize}
\item[(i)] if $|\Bx-\Bp| \leq C_1 \Ge$, then
\beq\label{near}
|\nabla u (\Bx)| \simeq \frac 1 {|\Bx-\Bp|^2},
\eeq

\item[(ii)] if $|\Bx-\Bp| \geq C_2 \Ge |\log \Ge|$, then
\beq\label{far}
|\nabla u (\Bx)| \simeq \frac 1 {\sqrt \Ge (\Ge + x_2^2)},
\eeq

\item[(iii)] if $C_1 \Ge < |\Bx-\Bp| < C_2 \Ge |\log \Ge|$, then
\beq\label{between}
|\nabla u (\Bx)| \lesssim \frac 1 {|\Bx-\Bp|^2} \exp \left(-A_* \frac {|\Bx-\Bp|}{ \Ge }\right)   + \Ge^{-3/2}.
\eeq
\end{itemize}
The constants involved in the relations $\simeq$ and $\lesssim$ above depend on $\Ba$ and $M$.
\end{thm}

We mention that \eqnref{near} and \eqnref{far} are estimates from below as well as from above, while \eqnref{between} is that from above and is a bridge between two estimates. \eqnref{near} shows that near the location of the emitter, the size of the field is of the same order as $|\nabla (\Ba \cdot \nabla \Ncal_\Bp(\Bx))|$. It is \eqnref{far} which exhibits field enhancement.   It is instructive to look into the estimate \eqnref{far} when $\Bp$ is located close to $\Bo$, the origin. If $|\Bp|\leq \sqrt{\Ge}$, then we see from \eqnref{far} that
\beq\label{far2}
|\nabla u (\Bx)| \simeq \frac 1 {\sqrt \Ge (\Ge + x_2^2)} \simeq \frac 1 {\sqrt \Ge |\Bx - \Bp|^2},
\eeq
provided that $ |\Bx| \geq 2\sqrt{\Ge}$.  Since the field $\nabla (\Ba \cdot \nabla \Ncal_\Bp(\Bx))$ excited by the emitter is of size $|\Bx-\Bp|^{-2}$, this inequality shows that the field is enhanced by the factor of $\Ge^{-1/2}$. It is quite interesting to observe that $\Ge^{-1/2}$ is the order of gradient blow-up when there is a smooth back-ground field, not emitter (see, e.g., \cite{AKL-MA-05, keller, Y}).

We have the following theorem when $\Ba=(0,1)$, which shows no enhancement of field.
\begin{thm}\label{thm_2nd}
Let $u$ be the solution to \eqref{main_equation}. If $\Ba =(0,1)$, then there exists a positive constant $A$ such that
\beq\label{urest}
|\nabla u (\Bx)| \lesssim \frac 1 {|\Bx - \Bp |^2} \exp \left(-A \frac {|\Bx-\Bp|}{\sqrt \Ge| \Bx- \Bp| + |p\Bx + (0,\Ge ) | }\right)
\eeq
for all $ \Bx = (x_1,x_2)\in B_{1/2} \setminus \overline{ D_1 \cup D_2}$, provided that $\Ge$ is sufficiently small and $|\Bp| <1/2$.
\end{thm}

In fact, \eqref{urest} shows not only that field enhancement does not occur in this case, but also that $|\nabla u|$ decays very fast. For example, if $|\Bp| \leq \sqrt \Ge$ and $|\Bx| > 2 \sqrt \epsilon$ for $\epsilon$ sufficiently small, then
$$
|\nabla u (\Bx)| \lesssim \frac 1 {|\Bx - \Bp|^2} \exp \left(-\frac A 3 \frac 1 {\sqrt \epsilon}\right)\leq \exp \left(-\frac A 4 \frac 1 {\sqrt \epsilon}\right).
$$

To prove these theorems, we decompose $u$ as
\beq\label{decomp_Q_r}
u = Q + r \quad \mbox {in } \Rbb^2 \setminus \overline{D_1 \cup D_2},
\eeq
where $Q$ and $r$ satisfy
\beq \label{eqn_Q}
\begin{cases}
\GD Q = 0  \quad&\mbox{in } \Rbb^2 \setminus \overline {(D_{1} \cup D_{2})}, \\
\ds Q = u~ (= c_j) \quad&\mbox{on }\p D_{j}, \  j=1,2,  \\
\ds \int_{\Rbb^2 \setminus \overline {(D_{1} \cup D_{2})} }  |\nabla Q|^2 d\Bx &< \infty,
\end{cases}
\eeq
and
\beq \label{eqn_r}
\begin{cases}
\GD r = {\bf a } \cdot \nabla  \Gd_{\Bp}  \quad&\mbox{in } \Rbb^2 \setminus \overline {(D_{1} \cup D_{2})}, \\
\ds r =0 \quad&\mbox{on }\p D_{j} , \  j=1,2,\\
\ds \int_{\Rbb^2 \setminus \overline {(D_{1} \cup D_{2})} }  |\nabla (r - \Ba \cdot \nabla
\Ncal_{\Bp}) |^2 d\Bx& < \infty.
\end{cases}
\eeq
We construct $Q$ and $r$ rather explicitly (see \eqref{Qr}), and estimate their gradients in the narrow region, say $B_{1/2} \setminus \overline{ D_1 \cup D_2}$.
One can see from the conditions on $\p D_j$ in \eqnref{eqn_Q} and \eqnref{eqn_r} that blow-up of $\nabla Q$ is caused by the potential difference on the closely located inclusions (and the emitter), while that of $\nabla r$ is caused solely by the existence of emitter. As the following proposition and Theorem \ref{thm_1st} show, if $\Ba \neq (0,1)$ and the emitter is located sufficiently close to the origin, then $\nabla Q$ dominates $\nabla u$ in almost all areas except a small neighborhood of the location of the emitter, which means that the field excited by the emitter is enhanced by the interaction of closely located inclusions.

\begin{prop}\label{main1}
The solution $u$ to \eqref{main_equation} admits the decomposition \eqref{decomp_Q_r} where $Q$ and $r$ satisfy \eqref{eqn_Q} and \eqref{eqn_r}, respectively. Moreover, there exists a positive constant $A$ such that the following inequalities hold for all sufficiently small $\Ge$, all $\Bp$ with $|\Bp| < 1/2$, all $ \Bx = (x_1,x_2)\in B_{1/2} \setminus \overline{ D_1 \cup D_2}$ and all unit vectors $\Ba=(a_1,a_2)$:
\beq\label{Qest}
\left| \nabla  Q (\Bx)\right| \simeq  |a_1|\frac {\sqrt \Ge}{( \Ge + p^2 )(\Ge + x_2^2)}
\eeq
and
\beq\label{rest}
|\nabla r (\Bx)| \lesssim \frac 1 {|\Bx - \Bp |^2} \exp \left(-A \frac {|\Bx-\Bp|}{\sqrt \Ge| \Bx- \Bp| + |p\Bx + (0,\Ge) | }\right).
\eeq
 The constants involved in the relations $\simeq$ and $\lesssim$ in the above can be chosen independently of $\Bp$ and $\Ba$ as well as $\Ge$.
\end{prop}

The estimate \eqnref{Qest} shows that that if $\Ba=(0,1)$, then $\nabla  Q (\Bx)=(0,0)$, thus \eqnref{urest} in Theorem \ref{thm_2nd} is nothing but \eqnref{rest}. The proof of Proposition \ref{main1} is presented in Section \ref{sec2} and that of Theorem \ref{thm_1st} is in the section that follows.

\section {Proof of Proposition \ref{main1}}\label{sec2}

We begin this section by recalling a special function. Let $R_j$ be the inversion with respect to the circle $\p D_j$ for $j=1,2$. Then the iterated inversions $R_1R_2$ and $R_2R_1$ have two fixed points $\Be_1 \in D_1$ and $\Be_2 \in D_2$. It is proved in \cite{KLY-JMPA-13} that
\beq\label{Bej}
\Be_j = \left((-1)^{j} \sqrt \Ge + O(\Ge), 0\right), \quad j=1,2.
\eeq
The following function was introduced in \cite{Y}:
\beq\label{expression_q}
q (\Bx) = \Ncal_{\Be_1}(\Bx) - \Ncal_{\Be_2}(\Bx) ,
\eeq
where $\Ncal_{\Be_j}$ is defined by \eqnref{Ncal}.
The function $q$ is harmonic in $\Rbb^2 \setminus \overline {(D_{1} \cup D_{2})}$ and attains a constant value on each boundary $\p D_j$ since $\p D_1$ and $\p D_2$ are Apollonius circles of $\Be_1$ and $\Be_2$. Thanks to symmetry of $D_1 \cup D_2$ with respect to the $x_2$-axis, we have
\beq\label{q_boundary_equality}
q |_{\p D_1 } = - q |_{\p D_2 }.
\eeq
Since $\GD \Ncal_\By =\Gd_\By$ and $q$ is constant on $\p D_j$, we see that
\beq\label{green}
\int_{\p D_{j} } v \,\p_{\nu} q \, ds = \int_{\p D_{j} } v \,\p_{\nu} q \, ds - \int_{\p D_{j} } \p_\nu v \, q \, ds= (- 1)^j v(\Be_j)
\eeq
for any harmonic function $v$ in $D_1 \cup D_2$. In particular, we have
\beq\label{thirdline}
\int_{\p D_{j} } \p_{\nu} q \, ds  =(-1)^j, \quad j=1,2.
\eeq
It is helpful to emphasize here that the normal vector $\nu$ is pointing inward to $D_1 \cup D_2$. Moreover, we have
\beq
q (\Bx) = O\left( |\Bx|^{-1}\right) \quad \mbox{as }|\Bx|\rightarrow\infty.
\eeq

We now briefly discuss on existence of the solution to \eqnref{main_equation}. Uniqueness of the solution can be easily shown using Green's theorem (or Hopf's lemma). Let $v$ be the unique solution to the problem
\beq\label{equation_of_v}
\begin{cases}
\GD v = 0  \quad&\mbox{in } \Rbb^2 \setminus \overline {(D_{1} \cup D_{2})}, \\
\ds v =  - \Ba \cdot \nabla\Ncal_{\Bp} \quad&\mbox{on }\p D_{1} \cup \p D_2,
\end{cases}
\eeq
satisfying
\beq\label{ltwo}
\int_{\Rbb^2 \setminus \overline {(D_{1} \cup D_{2})}} |\nabla v|^2 d\Bx  < \infty.
\eeq
The existence of $v$ can be shown using the Lax-Milgram theorem. Moreover, we infer from \eqnref{ltwo} that $|\nabla v(\Bx)| = O(|\Bx|^{-2})$ as $|\Bx| \to \infty$. Thus we have
\beq\label{meanzero}
\int_{\p D_1} \p_{\nu} v \, ds + \int_{\p D_2} \p_{\nu} v \, ds =0.
\eeq
Then one can see easily that the function $u$, defined by
\beq\label{uv}
u (\Bx)=  \Ba \cdot \nabla\Ncal_{\Bp}(\Bx) + v(\Bx) + c q (\Bx)- v_0,
\eeq
is the unique solution to \eqref{main_equation}. Here
\beq\label{cvzero}
c:= \int_{\p D_1} \p_{\nu} v \, ds \quad\mbox{and}\quad
v_0: = \lim_{|\Bx|\rightarrow \infty} v (\Bx).
\eeq

We see from \eqnref{uv} that
\beq \label{additional_Qr}
u|_{\p D_j} = cq|_{\p D_j} - v_0, \quad j=1,2.
\eeq
Thus we may take for the decomposition \eqnref{decomp_Q_r}
\beq\label{Qr}
Q=cq - v_0 \quad\mbox{and} \quad r= \Ba \cdot \nabla\Ncal_{\Bp} + v.
\eeq

\subsection{Proof of \eqnref{Qest}}

We see from \eqref{additional_Qr}  that
\beq\label{cid}
c= \frac { u |_{\p D_2 } - u |_{\p D_1}}{ q |_{\p D_2 } - q |_{\p D_1}}.
\eeq
We then adapt an argument in \cite{KLY-JMPA-13} to show that
\beq\label{podi}
u |_{\p D_2 } - u |_{\p D_1} = \left(\Ba \cdot \nabla\Ncal_{\Bp}\right)  (\Be_2) - \left(\Ba \cdot \nabla\Ncal_{\Bp}\right)  (\Be_1).
\eeq
In fact, by \eqnref{thirdline}, we have
$$
u |_{\p D_2 } - u |_{\p D_1} = \int_{\p D_1 \cup \p D_2} u \, \p_{\nu} q \, ds.
$$
We then have from \eqnref{uv} that
$$
u |_{\p D_2 } - u |_{\p D_1} = \int_{\p D_1 \cup \p D_2} \Ba \cdot \nabla\Ncal_{\Bp} \, \p_{\nu} q \,ds + \int_{\p D_1 \cup \p D_2} (v + c q -v_0) \, \p_{\nu} q \, ds.
$$
Using \eqnref{thirdline} and the definition of $c$ (in \eqnref{cvzero}), one can see that
$$
\int_{\p D_1 \cup \p D_2} q \p_{\nu} (v + c q - v_0)  \, ds=0.
$$
Then Green's theorem yields
$$
\int_{\p D_1 \cup \p D_2} (v + c q- v_0 ) \, \p_{\nu} q \, ds=0,
$$
and hence
$$
u |_{\p D_2 } - u |_{\p D_1} = \int_{\p D_1 \cup \p D_2} \Ba \cdot \nabla\Ncal_{\Bp} \, \p_{\nu} q \, ds .
$$
Now \eqnref{podi} follows by \eqnref{green}.

We see from \eqnref{podi} that
$$
u |_{\p D_2 } - u |_{\p D_1} = \frac{1}{2\pi} \frac{\Ba \cdot (\Be_2-\Bp)}{|\Be_2-\Bp|^2} - \frac{1}{2\pi} \frac{\Ba \cdot (\Be_1-\Bp)}{|\Be_1-\Bp|^2}.
$$
Thus one can show using \eqnref{Bej} that
\beq\label{udiffer}
\left|u |_{\p D_2 } - u |_{\p D_1}\right| \simeq |a_1| \frac{\sqrt {\Ge}} { {\Ge}+ p^2},
\eeq
if $\Ge$ is sufficiently small.

Using the explicit expression \eqnref{expression_q} of the function $q$ one can see that
\beq\label{qdiffer}
{ q |_{\p D_2 } - q |_{\p D_1}} \simeq \sqrt {\Ge}
\eeq
and
\beq\label{nablaq}
| \nabla q (\Bx)| \simeq \frac {\sqrt \Ge} { {\Ge} + x_2 ^2}
\eeq
for all $\Bx = (x_1,x_2) \in B_{1/2} \setminus (D_1 \cup D_2)$, if $\Ge$ is sufficiently small.

It now follows from \eqnref{cid}, \eqnref{udiffer} and \eqnref{qdiffer} that
$$
|c| \simeq \frac {|a_1|} { {\Ge} + x_2 ^2},
$$
which, together with \eqnref{nablaq}, leads us to \eqnref{Qest}.

\subsection{Proof of \eqnref{rest}}

Define the transformation $\Phi$ by
\beq\label{Phi}
\Phi (\By) =  \frac {\By-\Bp} {| {\By-\Bp}|^2 } +  \Bp,
\eeq
which enjoys the property that $\Phi(\Phi(\By))=\By$ for all $\By \neq \Bp$. Let
\beq
\GO_j := \Phi(D_j), \quad j=1,2.
\eeq
One can see that
\beq
\GO_{j} = \Ge_* D_j + \Bp_*  ,
\eeq
where
\beq\label{stars}
\Ge_* = \frac{1}{\Ge + p^2 + (\Ge^2 /4)} \quad\mbox{and}\quad \Bp_* = (0, p (1-\Ge_* )).
\eeq
Here and afterwards, $\Ge_* D_j$ denotes the dilation of $D_j$ by $\Ge_*$.

Recall that $r= \Ba \cdot \nabla\Ncal_{\Bp} + v$, and define $r_1$ by
\beq
r_1(\By):= r (\Phi(\By)).
\eeq
A straightforward computation shows that $\GD r_1=0$ in $\Rbb^2 \setminus  {(\GO_1 \cup \GO_2 \cup \{\Bp\})}$. Since
$$
\lim_{\By \to \Bp} r_1(\By) = \lim_{|\Bx| \to \infty} r(\Bx) = v_0,
$$
the point $\Bp$ is a removable singularity of $r_1$. Thus $r_1$ satisfies
\beq\label{r_varphi_1_eq}
\begin{cases}
\GD  r_1= 0 ~&\mbox{in } \Rbb^2 \setminus \overline {\GO_1\cup \GO_2}, \\
\ds  r_1 = 0 ~&\mbox{on }\p\GO_1  \cup \p \GO_2.
\end{cases}
\eeq
Moreover, \eqnref{ltwo} yields
$$
\int_{ \Rbb^2 \setminus \overline {\GO_1\cup \GO_2} } \left|\nabla \Big( r _1(\By)  -  \Ba \cdot (\nabla \Ncal_{\Bp})(\Phi(\By))\Big)  \right|^2 d\By < \infty .
$$
Since
$$
\Ba \cdot \nabla \Ncal_\Bp(\Phi(\By)) = \frac{1}{2\pi} \frac{\Ba \cdot (\Phi(\By)-\Bp)}{|\Phi(\By)-\Bp|^2} = \frac{1}{2\pi} \Ba \cdot (\By-\Bp),
$$
we have
\beq\label{r_varphi_2_eq}
\int_{ \Rbb^2 \setminus \overline {\GO_1\cup \GO_2} } \left| \nabla r_1(\By) - (1/2\pi) \Ba \right|^2 d\By < \infty.
\eeq

Define $r_2$ by
\beq
r_2 (\Bz)= r_1(\Ge_* \Bz + \Bp_*), \quad  \Bz \in \Rbb^2 \setminus \overline{(D_{1} \cup  D_{2})} . \label{relation_betw_By_Bz}
\eeq
Then \eqref{r_varphi_1_eq} and \eqref{r_varphi_2_eq} show that $r_2$ satisfies
\beq \label{equation_r_2}
\begin{cases}
\GD  r_2 = 0  \quad&\mbox{in } \Rbb^2 \setminus \overline{(D_{1} \cup  D_{2})}, \\
\ds r _2 = 0 \quad&\mbox{on }\p (D_{1} \cup  D_{2}),
\end{cases}
\eeq
with
\beq\label{condition_r_2}
\int_{ \Rbb^2 \setminus \overline {(D_{1} \cup  D_{2})}} \left|\nabla r_2 (\Bz) - (\Ge_*/2\pi )  \Ba \right|^2 d\Bz < \infty.
\eeq

Here we recall two known results. For a given harmonic function $h$ defined in $\Rbb^2$, let $u_h$ be the solution to
\beq \label{eqn_u_h}
\begin{cases}
\GD  u_h = 0  \quad&\mbox{in } \Rbb^2 \setminus \overline{(D_{1} \cup  D_{2})}, \\
\ds u_h = 0 \quad&\mbox{on }\p (D_{1} \cup  D_{2}),
\end{cases}
\eeq
with
\beq\label{condition_u_h}
\int_{ \Rbb^2 \setminus \overline {(D_{1} \cup  D_{2})}} \left| \nabla (u_h - h)  (\Bz) \right|^2 d\Bz < \infty.
\eeq
It is proved in \cite[Theorem 2.1]{KLY-JMPA-13} that for any $R>0$ there exists a constant $C_1$ such that
\beq\label{JMPA}
|\nabla u_h (\Bz) | \leq C_1
\eeq
for all $\Bz\in B_{R} \setminus \overline {(D_{1} \cup  D_{2})}$ and all sufficiently small $\Ge$. On the other hand, it is proved in \cite[Theorem 3]{KLY-MA-15} that there exist positive constants $A$, $\Gd$, and $C_2$ such that
\beq\label{MA}
|\nabla u_h(\Bz) | \leq C_2  \exp \left(- \frac A {\sqrt \Ge + |z_2 |}\right)
\eeq
for all $\Bz\in  B_{\Gd} \setminus \overline {D_{1} \cup  D_{2}} $ and all sufficiently small $\Ge$. Actually these results are obtained with the condition \eqnref{condition_u_h} replaced by
$$
u_h(\Bz) - h(\Bz) - a_h = O(|\Bz|^{-1}) \quad\mbox{as } |\Bz| \to \infty
$$ for some constant $a_h$.
However, the same proofs are valid even with the condition \eqnref{condition_u_h}.

From these results we obtain the following lemma.
\begin{lem}\label{lem:2circle}
There exists a positive constant $A$ independent of $\Ba$ such that
\beq\label{exp}
|\nabla r_2 (\Bz)| \lesssim \Ge_* \exp \left(-\frac A { {\sqrt \Ge} + |\Bz  | } \right)
\eeq
for all $\Bz\in  \Rbb^2 \setminus \overline{(D_{1} \cup  D_{2})}$ and all sufficiently small $\Ge$.
\end{lem}

\pf
Let $h (\Bz)= (1/\pi)\Ba \cdot \Bz $ so that $r_2 (\Bz)= \Ge_* u_h (\Bz)$ following notation in \eqnref{eqn_u_h}. By \eqnref{MA}, there are positive constants $A$,  $\Gd$, and $C_1$ such that
\beq\label{final_r_2_1}
|\nabla r_2(\Bz) | \leq C_1 {\Ge}_* \exp \left(- \frac A {\sqrt \Ge + |\Bz|}\right)
\eeq
for all $\Bz\in  B_{\Gd}  \setminus \overline {(D_{1} \cup  D_{2})}$.

Suppose that $R$ is large enough so that $B_R$ contains $\ol{D_1 \cup D_2}$. By \eqnref{JMPA},  there is $C_2>0$ such that  \beq\label{subfinal_r_2-1}
|\nabla r_2(\Bz) | \leq C_2  {\Ge}_*
\eeq
for all $\Bz\in  B_R  \setminus \overline {(D_{1} \cup  D_{2})}$. Thus we have
$$
\norm{ \nabla r_2- (\Ge_*/2\pi )  \Ba }_{L^{\infty}(\p B_R)} \leq \left(C_2  +  (1/{2\pi})\right)\Ge_* .
$$
Thanks to \eqref {condition_r_2}, we can apply the maximum principle on $\Rbb^2 \setminus B_R$ to obtain the following inequality for all $\Bz\in \Rbb^2 \setminus B_R$:
\begin{align}
|\nabla r_2(\Bz) | &\leq  | \nabla r_2(\Bz) - (\Ge_*/2\pi )  \Ba | +  | (\Ge_*/2\pi )  \Ba |\notag \\
&\leq  \norm{ \nabla r_2- (\Ge_*/2\pi )  \Ba }_{L^{\infty}(\p B_R)} + (\Ge_*/ 2\pi) \leq C_3 {\Ge}_* \label{subfinal_r_2-2}
\end{align}
with $C_3:= C_2 + 1/\pi$.

Let $A$ and $\Gd$ be constants appearing in \eqnref{final_r_2_1} and let $C_4 = \exp (A/\Gd)$. If $\Bz \in \Rbb^2 \setminus {B_\Gd}$, then
$$
1  \leq C_4 \exp (-  A/\Gd ) \leq C_4\exp \left(-  \frac A {|\Bz | }\right) \leq  C_4 \exp \left(- \frac A {\sqrt \Ge + |\Bz  | }\right).
$$
Thus we have from \eqref{subfinal_r_2-1} and \eqref{subfinal_r_2-2} that
$$
|\nabla r_2(\Bz) | \leq C_3 C_4 {\Ge_*}   \exp \left(- \frac A {\sqrt \Ge + |\Bz| }\right)
$$ for all $\Bz\in  \Rbb^2 \setminus (\overline {(D_{1} \cup  D_{2}) }   \cup B_{\Gd})$. This together with \eqref{final_r_2_1} yields \eqnref{exp}. Moreover, since \eqnref{exp} holds when $\Ba = (1,0)$ and $\Ba = (0,1)$, and $r_2$ depends on $\Ba$ linearly, we infer that $A$  satisfying \eqnref{exp}  can be chosen independently of $\Ba$. \qed

Recall that
\beq\label{ronetwo}
r(\Bx)=r_1(\By)=r_2(\Bz),
\eeq
where
\beq\label{BxByBz}
\Bx=\Phi(\By), \quad \By=\Ge_* \Bz+ \Bp_*
\eeq
with $\Ge_*$ and $\Bp_*$ defined in \eqnref{stars}. Since the complex conjugate of $\Phi$ is analytic, the Cauchy-Riemann equations hold:
$$
\frac {\p x_1}{\p y_1 }  = -\frac {\p x_2}{\p y_2} \quad\mbox{and}\quad \frac {\p x_1}{\p y_2 }  = \frac {\p x_2}{\p y_1}.
$$
Thus, we have
$$
\left(\frac {\p x_1}{\p y_1} , \frac {\p x_2}{\p y_1}\right) \cdot \left(\frac {\p x_1}{\p y_2} , \frac {\p x_2}{\p y_2}\right) = 0 ,
$$
and
\begin{align}
& \sqrt {\left|\left(\frac {\p x_1}{\p y_1} , \frac {\p x_2}{\p y_1}\right) \right|^2 +   \left|\left(\frac {\p x_1}{\p y_2} , \frac {\p x_2}{\p y_2}\right) \right|^2 } \nonumber\\
& = \sqrt 2 \left|\left(\frac {\p x_1}{\p y_1} , \frac {\p x_1}{\p y_2}\right) \right|
= \sqrt 2 |\By - \Bp|^{-2}
=  \sqrt 2 |\Bx - \Bp|^2. \label{sqrt}
\end{align}
We see from \eqnref{ronetwo} and \eqnref{BxByBz} that
\begin{align*}
|\nabla r_2 (\Bz)| & =  \left|\left(\frac{\p r_1} {\p y_1} \left(\By\right) \frac {\p y_1}{\p z_1} , \frac{\p r_1} {\p y_2}\left(\By\right) \frac {\p y_2}{\p z_2}\right) \right| \\
& =  {\Ge}_* \left|\nabla \left(r \left(\Phi(\By)\right)\right)\right|\\
& =  {\Ge}_* \left|\nabla r (\Bx)\right|
\sqrt {\left|\left(\frac {\p x_1}{\p y_1} , \frac {\p x_2}{\p y_1}\right) \right|^2 +
\left|\left(\frac {\p x_1}{\p y_2} , \frac {\p x_2}{\p y_2}\right) \right|^2 }.
\end{align*}
It then follows from \eqnref{sqrt} that
\beq
|\nabla r_2 (\Bz)| = \sqrt 2   {\Ge}_* |\Bx - \Bp|^2 \left|\nabla r (\Bx)\right|.
\eeq
Then \eqnref{exp} yields
$$
\left|\nabla r (\Bx)\right| \lesssim \frac 1 {|\Bx-\Bp|^2} \exp \left(-A \frac 1 { {\sqrt \Ge} + |\Bz  | } \right).
$$
Note that
\begin{align*}
\left |\Bz \right| & = \left| \Ge_*^{-1} (\By - \Bp_*)   \right|  = \left| \Ge_*^{-1} (\By - \Bp) +(0,p )  \right| \\& = \left| \Ge_*^{-1}  \frac { \Bx - \Bp } { |\Bx - \Bp|^2 }  + (0,p )  \right| \\&= |\Bx - \Bp|^{-1} |p(\Bx - \Bp) + (0,\Ge_*^{-1} ) | =  |\Bx - \Bp|^{-1} |p\Bx + (0, \Ge+ \Ge^2/4) |.
\end{align*}
Thus
\beq\label{r_gradient_estimate_semi_final}
|\nabla r (\Bx)| \lesssim \frac 1 {|\Bx - \Bp |^2} \exp \left(-A \frac {|\Bx-\Bp|}{\sqrt \Ge| \Bx- \Bp| + |p\Bx + (0,\Ge + \Ge^2/4) | }\right).
\eeq
Now \eqnref{rest} follows from the following lemma and the proof is complete.

\begin{lem}
It holds that
\beq
\sqrt \Ge| \Bx- \Bp| +  |p\Bx + (0,\Ge + \Ge^2/4) | \simeq \sqrt \Ge| \Bx- \Bp| +  |p\Bx + (0,\Ge) |
\eeq
for all $\Bx \in \Rbb^2$, all $\Ge < 1/2$ and all $\Bp=(0,p)$.
\end{lem}
\pf
Since
$$
|x_1| + |x_2| \simeq |(x_1,x_2)|
$$
for all $\Bx  = (x_1,x_2)\in \Rbb^2$ and $\Bp$ is of the form $(0,p)$, we have
\begin{align*}
\sqrt \Ge| \Bx- \Bp| +  |p\Bx + (0,\Ge + \Ge^2/4) |  &\simeq (\sqrt \Ge + |p |) |x_1| + \left(\sqrt \Ge |x_2-p| + |px_2+ \Ge + \Ge^2/4|\right),\\
\sqrt \Ge| \Bx- \Bp| +  |p\Bx + (0,\Ge) |  &\simeq (\sqrt \Ge + |p |) |x_1| + (\sqrt \Ge |x_2-p| + |px_2+ \Ge |).
\end{align*}
Thus it suffices to show that
\beq\label{lemma_equv_norm:1}
\sqrt \Ge| x_2- p| + |px_2 + \Ge + \Ge^2/4 | \simeq \sqrt \Ge | x_2- p| + |px_2 + \Ge |.
\eeq

We now consider two cases separately: when $\sqrt \Ge |x_2 -p |\geq \Ge^{2}$, and when $\sqrt \Ge |x_2 -p |< \Ge^{2}$.
The triangular inequality yields the following inequalities:
$$
\sqrt \Ge | x_2- p| - \Ge^2 /4 + |px_2 +\Ge | \le \sqrt \Ge| x_2- p| + |px_2 + \Ge + \Ge^2/4 | \le \sqrt \Ge| x_2- p| + \Ge^2 /4 + |px_2 +\Ge |.
$$
Moreover, if $\sqrt \Ge| x_2- p|  \geq  \Ge^2$, we have
\begin{align*}
\sqrt \Ge| x_2- p| +(\Ge^2 /4)+ |px_2 +\Ge |  &\leq 2 \left(\sqrt \Ge| x_2- p| + |px_2 +\Ge |\right), \\
\sqrt \Ge| x_2- p| - (\Ge^2 /4)+ |px_2 +\Ge |  &\geq (3/4) \left(\sqrt \Ge| x_2- p| + |px_2 +\Ge |\right).
\end{align*}
Thus, \eqref{lemma_equv_norm:1} holds in the first case, namely, when $\sqrt \Ge |x_2 -p |\geq \Ge^{2}$.

In the second case, we prove that
\beq\label{lemma_equv_norm:1:2nd_case}
\left|px_2 + \Ge + \Ge^2/4 \right| \simeq  \left|px_2 + \Ge  \right|,
\eeq
which clearly implies the desired estimate \eqref{lemma_equv_norm:1}.

To prove \eqnref{lemma_equv_norm:1:2nd_case}, we start with inequalities
\beq\label{lemma_equv_norm:3}
p^2 + \Ge - \left( |p(x_2 - p)| +  \Ge^2/4 \right) \leq \left|px_2 + \Ge + \Ge^2/4 \right| \leq  p^2 + \Ge + \left(|p(x_2 - p)| + \Ge^2/4 \right),
\eeq
which are consequences of the triangular inequality. Since $\sqrt{\Ge} |x_2 -p |< \Ge^{2}$, we have
$$
|p(x_2 - p)| + \Ge^2/4 \leq  (p^2 + |x_2 -p |^2)/2 + \Ge^2/4 \leq (p^2 + \Ge^3)/2 + \Ge^2 /4.
$$
So, if $\Ge< 1/2$, then
$$
|p(x_2 - p)| +( \Ge^2/4) \leq 1/2 \left(p^2 + \Ge  \right) .
$$
Thus, \eqnref{lemma_equv_norm:3} yields
$$
1/2 \left(p^2 + \Ge  \right) \leq \left|px_2 + \Ge + \Ge^2/4 \right| \leq  3/2 \left(p^2 + \Ge  \right).
$$
Since $|x_2 -p |< \Ge^{3/2}$, one can easily that
$$
p^2 + \Ge \simeq  \left|px_2 + \Ge  \right|.
$$
Thus \eqref{lemma_equv_norm:1} follows in the second case, and the proof is complete. \qed
\section {Proof of Theorem \ref{thm_1st}}\label{sec_cor_1st}

We first recall that $\Bp$ satisfies the following condition for some $M$:
\beq\label{Bp}
|\Bp| < M \sqrt{\Ge}.
\eeq

If $|\Bx -\Bp| \le \Ge/4$, then $|\Bx| < (M + 1/4) \sqrt {\Ge}$ for all sufficiently small $\Ge$, and hence we infer from \eqnref{Qest} that
\beq\label{Qone}
|\nabla Q (\Bx)| \simeq \Ge^{-3/2} .
\eeq

Let $u_0(\Bx):= \Ba \cdot \nabla \Ncal_\Bp(\Bx)$ for ease of notation, and write $u$ as
$$
u (\Bx)= Q (\Bx) + \left( r(\Bx) - u_0(\Bx) \right) + u_0(\Bx).
$$
We see from \eqnref{rest} that
\beq\label{rone}
|\nabla r(\Bx)| \lesssim \Ge^{-2}
\eeq
for all $\Bx$ satisfying $|\Bx-\Bp|=\Ge/4$. Note that
$$
\nabla u_0 = \frac{1}{2\pi |\Bx-\Bp|^2} \left[ \Ba - \frac{2 \Ba\cdot (\Bx-\Bp) (\Bx-\Bp)}{|\Bx-\Bp|^2} \right].
$$
Moreover, one can easily see that
$$
\left| \Ba - \frac{2 \Ba\cdot (\Bx-\Bp) (\Bx-\Bp)}{|\Bx-\Bp|^2} \right| =|\Ba|=1
$$
for all $\Bx \neq \Bp$, and hence
\beq\label{uone}
|\nabla u_0 (\Bx)| = \frac 1 {2\pi |\Bx- \Bp|^2}.
\eeq
In particular, if $|\Bx-\Bp| = \Ge/4$, then
\beq\label{utwo}
|\nabla u_0(\Bx)| = \frac{8}{\pi}\Ge^{-2}.
\eeq
Since $\GD u_0= \Ba \cdot \nabla \Gd_{\Bp}$, we have $\GD(r-u_0)(\Bx)=0$ if $|\Bx-\Bp| \le  \Ge/4$. Moreover, we have from \eqnref{rone} and \eqnref{utwo} that
$$
\left|\nabla (r-u_0)(\Bx) \right| \le \left|\nabla r(\Bx) \right| + \left|\nabla u_0(\Bx) \right| \lesssim \Ge^{-2}
$$
if $|\Bx-\Bp| =  \Ge/4$.
Thus the maximum principle yields
\beq\label{ruone}
\left|\nabla (r-u_0)(\Bx) \right| \lesssim \Ge^{-2}
\eeq
for all $\Bx$ such that $|\Bx-\Bp| \le  \Ge/4$.

Now we infer using \eqnref{Qone}, \eqnref{uone} and \eqnref{ruone} that there exists a constant $C_1 \le 1/4$ such that
$$
|\nabla u (\Bx)| \simeq \frac 1 {|\Bx- \Bp|^2}
$$
for all $\Bx$ with $|\Bx-\Bp| < C_1\Ge$. This proves \eqnref{near}.

If $ |\Bx -\Bp | \geq 2 M \sqrt \Ge $, then we have $|\Bx| \geq M \sqrt \Ge$ thanks to \eqnref{Bp}, and hence
$$
2 |\Bx-\Bp| \ge |\Bx-\Bp| + |\Bp|\ge  |\Bx|.
$$
Since $a_1\neq 0 $, it follows from \eqnref{Qest} and \eqnref{Bp} that
\beq\label{Qthree}
|\nabla Q (\Bx)| \simeq \frac {1}{ \sqrt\Ge (\Ge +  x_2^2)} ,
\eeq
and from \eqnref{rest} that
$$
|\nabla r (\Bx)| \lesssim \frac {1}{|\Bx- \Bp|^2} \lesssim  \frac {1}{ |\Bx|^2} \leq \frac 2 {M^2 \Ge + x_2^2}.
$$
Thus, if $|\Bx -\Bp | \geq 2 M \sqrt \Ge$, we have
\beq\label{uthree}
|\nabla u (\Bx)| \simeq \frac {1}{ \sqrt\Ge (\Ge +  x_2^2)}
\eeq
for all sufficiently small $\Ge$.

Now suppose that $0 <|\Bx -\Bp | \leq 2M \sqrt \Ge$. Then $|x_2|\le |\Bx| \le 3M \sqrt{\Ge}$. Thus we have from \eqnref{Qthree} that
\beq\label{Qfour}
|\nabla Q (\Bx)| \simeq \frac {1}{ \Ge \sqrt\Ge }  .
\eeq
Moreover, one can see that
$$
\sqrt {\Ge }|\Bx - \Bp| +|p \Bx + (0,\Ge )| \le (2M+ 3M^2 + 1 )\Ge .
$$
Thus it follows from \eqnref{rest} that
\beq\label{rtwo}
|\nabla r (\Bx)| \lesssim \frac {1}{ |\Bx-\Bp |^2} \exp \left(-A_* \frac {|\Bx - \Bp|}{\Ge  }\right) ,
\eeq
where $A_*$ is the constant defined by
$$
A_* = \frac A {2M+ 3M^2 + 1}.
$$

If further $C_2 \Ge |\log \Ge | \leq |\Bx-\Bp| \leq 2 M \sqrt \Ge$ with $C_2 =  \frac 1 {2A_*}$, then it follows from \eqnref{rtwo} that
$$
|\nabla r (\Bx)| \lesssim   \frac 1 {|\Bx-\Bp|^2} \exp (- (\log \Ge)/2) \leq \frac 1 {\Ge \sqrt {\Ge} |\log \Ge|^2}.
$$
We then see from \eqnref{Qfour} that
\beq\label{ufour}
|\nabla u (\Bx)| \simeq |\nabla Q(\Bx)| \simeq \frac 1 {\Ge \sqrt \Ge} \simeq \frac 1 { \sqrt \Ge (\Ge+ x_2^2)}
\eeq
for all sufficiently small $\Ge$. Now \eqnref{far} follows from \eqnref{uthree} and \eqnref{ufour}.

The inequality \eqref{between} is an immediate consequence of \eqnref{Qfour} and \eqnref{rtwo}. This completes the proof.  \qed

\section*{Conclusion}

We study enhancement, in the presence of closely located circular inclusions, of the field excited by an emitter and derive precise estimates quantifying such enhancement. Estimates show that the field is enhanced at points away from the emitter and it is due to strong interaction between inclusions. They also show that the magnitude of enhancement is of the factor $\Ge^{-1/2}$ where $\Ge$ is the distance between two inclusions. This factor is the same as that for the field enhancement in the case that a smooth back-ground field, not an emitter, is present. In the companion paper \cite{KY2} the field enhancement is considered in the presence of a bow-tie structure, and it is showed with  precise estimates that the field is enhanced near vertices due to corner singularities.

It is likely that the field is enhanced by the factor of $\Ge^{-1/2}$ for general strictly convex inclusions with smooth boundaries in two dimensions. It would be quite interesting to clarify this. Field enhancement in presence of an emitter and spherical inclusions in three dimensions is also studied to confirm that the factor of the enhancement is $(\Ge |\log \Ge|)^{-1}$. This factor is in accordance with results obtained in \cite{BLY-ARMA-09, LY}. This result will be reported in a forthcoming paper.

 \end{document}